\documentclass[a4paper]{article}     
\usepackage{theorem}                 
\usepackage{amssymb}
\usepackage{epsfig}
\newtheorem{theorem}{Theorem}

\newtheorem{proposition}{Proposition}

\newcommand{\CC}{{\mathbb{C}}}
\newcommand{\HH}{{\mathbb{H}}}
\newcommand{\PP}{{\mathbb{P}}}

\newcommand{\ZZ}{{\mathbb{Z}}}

\newcommand{\NN}{{\mathbb{N}}}

\begin{document}

\title{Poincar\'{e} series and monodromy of a two-dimensional quasihomogeneous
hypersurface singularity}
\author{Wolfgang Ebeling\\
\parbox{9cm}{\small
 \begin{center} Institut f\"{u}r Mathematik, Universit\"{a}t Hannover, \\ 
        Postfach 6009, D-30060 Hannover, Germany \\
        E-mail: ebeling@math.uni-hannover.de
\end{center}}
}
\date{}

\maketitle

\begin{abstract}
A relation is proved between the Poincar\'{e} series of the coordinate algebra of a
two-dimensional quasihomogeneous isolated hypersurface singularity and the characteristic polynomial
of its monodromy operator. For a Kleinian singularity not of type $A_{2n}$, this amounts to the
statement that the Poincar\'{e} series is the quotient of the characteristic polynomial of the
Coxeter element by the characteristic polynomial of the affine Coxeter element of the corresponding
root system. We show that this result also follows from the McKay correspondence. 
\end{abstract}

\section*{Introduction}

S.~M.~Gusein-Zade, F.~Delgado, and A.~Campillo \cite{GZDC99} have
shown that for an irreducible plane curve singularity the Poincar\'{e} series
of the ring of functions on the curve coincides with the zeta function of its monodromy
transformation.

In this paper we show that there is also a relation between the Poincar\'{e} series of the coordinate
algebra of a two-dimensional quasihomogeneous isolated hypersurface singularity and the characteristic
polynomial of its monodromy operator.

Let $(X,x)$ be a normal surface singularity  with good
$\CC^\ast$-action. The coordinate algebra $A$ is a graded algebra. We consider the
Poincar\'{e} series $p_A(t)$ of $A$. Let $\{ g; b; (\alpha_1, \beta_1), \ldots , (\alpha_r, \beta_r)
\}$ be the orbit invariants of $(X,x)$. 
We define
\begin{eqnarray*}
\psi_A(t) & := & (1-t)^{2-r}\prod_{i=1}^r (1-t^{\alpha_i}), \\
\phi_A(t) & := & p_A(t)\psi_A(t).
\end{eqnarray*}
Let $(X,x)$ be a
hypersurface singularity. Then $\phi_A(t)$ is a product of cyclotomic polynomials.

K.~Saito \cite{Saito94, Saito98} has introduced a duality between polynomials which are
products of cyclotomic polynomials. He has shown that V.~I.~Arnold's strange duality
between the 14 exceptional unimodal hypersurface singularities is related to such a
duality between the characteristic polynomials of the monodromy operators of the
singularities. It is now well-known
that Arnold's strange duality is related to the mirror symmetry of $K3$ surfaces (see
e.g.\
\cite{Dolgachev96}).

The main results of the paper are the following.
We show that the dual (in
Saito's sense) of the rational function
$\tilde{\phi}_A(t):=\phi_A(t)/(1-t)^{2g}$ is the characteristic polynomial of the monodromy
operator of $(X,x)$ (Theorem~\ref{phi_M3}). Similar results can be proved for isolated complete
intersection singularities (abbreviated ICIS in the sequel) of certain types (see
Theorem~\ref{phi_M4a} and \ref{phi_M4b}). 

If $(X,x)$ is a Kleinian singularity not of type $A_{2n}$, then $\psi_A(t)$ is the characteristic
polynomial of the affine Coxeter element of the corresponding root system and the above result implies
that $\phi_A(t)$ is the characteristic polynomial of the Coxeter element. Hence the Poincar\'{e}
series of a polyhedral group which is not a cyclic group of odd order is the quotient of these two
polynomials. We derive this result also directly from the McKay correspondence using ideas of the
paper 
\cite{Kostant85}. There are various formulas for Poincar\'{e} series of binary polyhedral groups in 
\cite{Knorrer85, Kostant85, Springer87}, but this relation seems to be new. 

The author thanks K.~Hulek for useful discussions. He is grateful to E.~Bries\-korn for
pointing out an error in an earlier version of the paper.

\section{Main results} \label{QSS}

Let $(X,x)$ be a normal surface singularity with a good $\CC^\ast$-action. So $X$ is a
normal two-dimensional affine algebraic variety over
$\CC$ which is smooth outside its {\em vertex} $x$. Its coordinate ring $A$ has the
structure of a graded $\CC$-algebra 
$A = \bigoplus_{k=0}^\infty A_k$, $A_0=\CC$, and $x$ is defined by the maximal ideal
$\frak{m}= \bigoplus_{k=1}^\infty A_k$. 

According to I.~Dolgachev \cite{Dolgachev75}, there exist a simply connected Riemann
surface $\cal D$, 
a discrete cocompact subgroup $\Gamma$ of $\mbox{Aut}({\cal D})$ and a line bundle
$\cal L$ on $\cal D$ to which the action of $\Gamma$ lifts such that 
$$A_k = H^0({\cal D}, {\cal L}^k)^\Gamma.$$

Let $Z:={\cal D}/\Gamma$. By \cite[Theorem~5.1]{Pinkham77a} (see also
\cite[Theorem~5.4.1]{Wagreich81}), there exist a divisor
$D_0$ on $Z$, $p_1, \ldots , p_r \in Z$, and integers $\alpha_i$, $\beta_i$
with $0<\beta_i < \alpha_i$ and $(\alpha_i,\beta_i)=1$ for 
$i=1, \ldots, r$ such that 
$$A_k =L \left( kD_0 + \sum_{i=1}^r \left[ k \frac{\alpha_i - \beta_i}{\alpha_i} \right]
p_i \right).$$
Here $[x]$ denotes the largest integer $\leq x$, and $L(D)$ for a divisor $D$ on $Z$
denotes the linear space of meromorphic functions $f$ on $Z$ such that $(f) \geq -D$.
We number the points $p_i$ so that $\alpha_1
\leq \alpha_2 \leq \ldots \leq \alpha_r$. Let $g$ be the genus of $Z$ and define $b :=
{\rm degree}\, D_0 +r$. Then $\{ g; b; (\alpha_1, \beta_1), \ldots , (\alpha_r, \beta_r)
\}$ are called the {\em orbit invariants} of $(X,x)$, cf.\ e.g.\ \cite{Wagreich83}.
Define
$\mbox{vdeg}({\cal L}):= -b + \sum_{i=1}^r \frac{\beta_i}{\alpha_i}$.

Now assume that $(X,x)$ is Gorenstein. By \cite{Dolgachev83}, there exists an integer
$R$ such that ${\cal L}^{-R}$ and the tangent bundle $T_{\cal D}$ of $\cal D$ are
isomorphic as $\Gamma$-bundles and
\begin{eqnarray*}
R \cdot \mbox{vdeg}({\cal L}) & = & 2 -2g - r + \sum_{i=1}^r \frac{1}{\alpha_i}, \\
R\beta_i & \equiv & 1 \ \mbox{mod} \, \alpha_i, \quad i=1, \ldots , r.
\end{eqnarray*}
Following \cite[3.3.15]{Dolgachev??} we call $R$ the {\em exponent} of $(X,x)$. 
Since $b$ and the $\beta_i$ are
determined by the $\alpha_i$ and the number $R$, we write the
orbit invariants also as $g; \alpha_1, \ldots, \alpha_r$.

The Gorenstein surface 
singularities with good $\CC^\ast$-action fall into three classes \cite{Dolgachev83,
Wagreich83}.

(1) ${\cal D} = \PP^1(\CC)$: Then $R=-2$ or $R=-1$, $g=0$, and $(X,x)$ is a Kleinian
singularity. 

(2) ${\cal D} = \CC$: Then $R=0$, $r=0$, and $g=1$. Hence $(X,x)$ is a
simply elliptic singularity \cite{Saito74}.

(3) ${\cal D} = \HH$, the upper half plane: The remaining Gorenstein surface
singularities with good $\CC^\ast$-action belong to this class. We have $R \geq 1$.  

We consider the {\em
Poincar\'{e} series} of the algebra $A$ 
$$ p_A(t)= \sum_{k=0}^\infty a_k t^k $$
where $a_k=\dim A_k$. It is well known that $p_A(t)$ is a rational function and the
the order of the pole of $p_A(t)$ at $t=1$ is equal to the dimension of $A$, hence equal
to 2. Moreover, $p_A(t)$ has simple poles at the $\alpha_i$-th roots of unity different
from 1. 

We define
\begin{eqnarray*}
\psi_A(t) & := & (1-t)^{2-r}\prod_{i=1}^r (1-t^{\alpha_i}), \\
\phi_A(t) & := & p_A(t)\psi_A(t), \\
\tilde{\phi}_A(t) & := & \frac{\phi_A(t)}{(1-t)^{2g}}.
\end{eqnarray*}
Then $\phi_A(t)$ is a polynomial. 

Now let $(X,x)$ be an ICIS
with weights $q_1, \ldots, q_n$ and degrees $d_1, \ldots , d_{n-2}$. Then its
Poincar\'{e} series is given by (see e.g.\ \cite[Proposition~(2.2.2)]{Wagreich83})
$$p_A(t) = \frac{\prod_{i=1}^{n-2} (1-t^{d_i})}{\prod_{j=1}^n (1-t^{q_j})}.$$
Hence $p_A(t)$, $\psi_A(t)$, $\phi_A(t)$, and $\tilde{\phi}_A(t)$ are rational functions of the form
$$\phi(t)= \prod_{m|h} (1 - t^m)^{\chi_m} \quad \mbox{for} \
\chi_m \in \ZZ \mbox{ and for some } h \in \NN.$$

Given a rational function
$$\phi(t)=\prod_{m|h} (1 - t^m)^{\chi_m},$$
K.~Saito \cite{Saito94} has defined a dual rational function 
$$\phi^\ast(t) = \prod_{k | h} (1-t^k)^{-\chi_{h/k}}.$$

We are now able to state the main results of the paper. 

\begin{theorem} \label{phi_M3}
Let $(X,x)$ be a quasihomogeneous hypersurface singularity in $\CC^3$. 
Then $\tilde{\phi}_A^\ast(t)$ is the characteristic polynomial of the classical
monodromy operator of $(X,x)$.
\end{theorem}

The proof of Theorem~\ref{phi_M3} will be given in Section~\ref{MO}.

Let $(X,x)$ be a Kleinian singularity. Then $g=0$, and therefore $\tilde{\phi}_A(t)=\phi_A(t)$. It is
well known that the Kleinian singularities correspond to root systems of type $A_l$, $D_l$, $E_6$,
$E_7$, or $E_8$. The classical monodromy operator of $(X,x)$ is the Coxeter element of the
corresponding root system.   
The polynomial $\psi_A(t)$ is the characteristic polynomial of the affine Coxeter element of the
corresponding root system (see \cite[p.~591]{Steinberg85} or \cite[6.2]{Springer87}, if $(X,x)$ is
not of type $A_l$). In the case $A_{2n}$, we have
$$\phi_A(t)= \frac{1-t^{4n+2}}{1-t^2}.$$
In this case, $\phi_A^\ast(t) \neq \phi_A(t)$. In all other cases one can verify that
$\phi_A^\ast(t) = \phi_A(t)$ (cf.\ Table~\ref{table4.1}). Therefore we obtain from
Theorem~\ref{phi_M3}: 

\begin{theorem} \label{McKay}
Let $(X,x)$ be a Kleinian singularity not of type $A_{2n}$. Then  $\phi_A(t)$ and $\psi_A(t)$ are 
the characteristic polynomials of the Coxeter element and 
the affine Coxeter element respectively of the corresponding root system. Hence $p_A(t)$ is the
quotient of these polynomials.
\end{theorem}

In Section~\ref{KS} we shall give a direct proof of Theorem~\ref{McKay} using the McKay
correspondence.

\section{Poincar\'{e} series and monodromy}
\label{MO}

In this section we shall prove Theorem~\ref{phi_M3}.

Let $(X,x)$ be a Gorenstein surface singularity with a good $\CC^\ast$-action.
The residue of the Poincar\'{e} series $p_A(t)$ at a primitive $\alpha_i$-th root of unity can be
computed as follows.

\begin{proposition} \label{prop:res}
Let $(X,x)$ be Gorenstein and $R$ be the exponent of $(X,x)$.
Let $\xi_i= \exp (2\pi\sqrt{-1}/\alpha_i)$. Then the residue of $p_A(t)$ at $t=\xi_i$ is
equal to 
$$\sum_{\alpha_i|\alpha_j}\frac{\xi_iÊ\cdot \xi_i^R}{\alpha_j(1-\xi_i^R)}.$$
\end{proposition}

From Proposition~\ref{prop:res} we can derive the following proposition generalizing \cite[Proposition
(2.8)]{Wagreich80}).
For integers $a_1, \ldots , a_r$ we denote by $\langle a_1, \ldots ,a_r
\rangle$  their least common multiple and by
$(a_1, \ldots , a_r)$ their greatest common divisor.

\begin{proposition} \label{Wagreich3}
Let $A$ be the coordinate algebra of a Gorenstein surface singularity with good
$\CC^\ast$-action. Suppose that $A$ is generated by 3 elements of weights $q_1$, $q_2$,
$q_3$. For each
$i=1, \ldots , r$ let $\xi_i= \exp (2\pi\sqrt{-1}/\alpha_i)$.
Then the Poincar\'{e} series is given by
$$p_A(t) = \frac{(1-t^d)}{(1-t^{q_1})(1-t^{q_2})(1-t^{q_3})}$$
if and only if the following conditions hold:

{\rm (a)} $(q_1,q_2,q_3) =1$.

{\rm (b)} $2g-2+r - \sum_{i=1}^r \frac{1}{\alpha_i}= 
\frac{Rd}{q_1q_2q_3}$.

{\rm (c)} For each $i$
$$\sum_{\alpha_i|\alpha_j}\frac{\xi_iÊ\cdot \xi_i^R}{\alpha_j(1-\xi_i^R)} = \left\{
\begin{array}{ll}
    \frac{-\xi_i(1-\xi_i^d)}{q_{i_1}(1-\xi_i^{q_{i_2}})(1-\xi_i^{q_{i_3}})} & \mbox{if
} \alpha_i|q_{i_1} \mbox{ and } \alpha_i \! \not| q_{i_2}, q_{i_3}, \\
    \frac{-d\xi_i}{q_{i_1}q_{i_2}(1-\xi_i^{q_{i_3}})} & \mbox{if } \alpha_i|q_{i_1},
\alpha_i|q_{i_2}, \alpha_i \! \not| q_{i_3}. \end{array} \right.$$

{\rm (d)} For all $i$ and $j$ so that $i\neq j$, $(q_i,q_j)|d$.

{\rm (e)} $d=q_1+q_2+q_3+R$.
\end{proposition}

Let $(X,x)$ be an isolated hypersurface singularity in $\CC^3$ given by a quasihomogeneous
equation $f(z_1,z_2,z_3)=0$ of degree $d$ and weights $q_1$, $q_2$, $q_3$. Then the  
Poincar\'{e} series is given by
$$p_A(t) = \frac{(1-t^d)}{(1-t^{q_1})(1-t^{q_2})(1-t^{q_3})}.$$
From Proposition~\ref{Wagreich3}
one can easily derive the following proposition which was proven by P.~Orlik and
P.~Wagreich 
\cite[3.6~Proposition~1]{OW77} using another method.

\begin{proposition}[Orlik,Wagreich] \label{OW}
Let $(X,x)$ be a quasihomogeneous isolated hypersurface singularity. Let $w_i= d/q_i = u_i/v_i$
where $(u_i,v_i)=1$ und $u_i, v_i \geq 1$. Assume $1
\leq v_1 \leq v_2 \leq v_3$. Then the table below indicates the number of orbit invariants $\alpha$ of
each type:
\begin{center}
\begin{tabular}{ccccccc}
$\alpha=$ & $(q_2,q_3)$ & $(q_1,q_3)$ & $(q_1,q_2)$ & $q_3$ & $q_2$ & $q_1$ \\
$1=v_1=v_2=v_3$ & $\frac{d}{\langle q_2,q_3 \rangle}$ & $\frac{d}{\langle q_1,q_3 \rangle}$ &
$\frac{d}{\langle q_1,q_2 \rangle}$ & & &
\\
$1=v_1=v_2<v_3$ & $\frac{d-q_2}{\langle q_2,q_3 \rangle}$ & $\frac{d-q_1}{\langle q_1,q_3 \rangle}$ &
$\frac{d}{\langle q_1,q_2 \rangle}$ & $1$ & &
\\ $1=v_1 < v_2 \leq v_3$ & $\frac{d-q_2-q_3}{\langle q_2,q_3 \rangle}$ & $\frac{d-q_1}{\langle
q_1,q_3 \rangle}$ &
$\frac{d-q_1}{\langle q_1,q_2 \rangle}$ & $1$ & $1$ & \\
$1 < v_1 \leq v_2 \leq v_3$ & $\frac{d-q_2-q_3}{\langle q_2,q_3 \rangle}$ & $\frac{d-q_1-q_3}{\langle
q_1,q_3 \rangle}$ &
$\frac{d-q_1-q_2}{\langle q_1,q_2 \rangle}$ & $1$ & $1$ & $1$
\end{tabular}
\end{center} The blank entries are zero if $q_i$ does not divide $q_j$ for $j \neq i$. If $q_i|q_j$,
then
$(q_i,q_j)=q_i$ and we list those orbit invariants under the column headed $(q_i,q_j)$.
\end{proposition} 

On the other hand, we consider the characteristic polynomial of the monodromy operator of $(X,x)$.
Let ${\cal X}_t$ be a Milnor fibre of the singularity $(X,x)$ and denote by
$M:=H_2({\cal X}_t, \ZZ)$ the corresponding Milnor lattice. 
Let $c: M \to M$ be the classical monodromy operator of the singularity $(X,x)$. It is
well-known that
$c$ is quasi-unipotent and therefore the eigenvalues of
$c$ are roots of unity. We write the characteristic polynomial $\phi_M(t)=\det
(tI-c)$ of $c$ as
$$\phi_M(t)=\prod_{m|h} (t^m-1)^{\chi_m} \quad \mbox{for} \
\chi_m \in \ZZ \mbox{ and for some } h \in \NN.$$

The characteristic polynomial
$\phi_M(t)$ can be computed as follows \cite{MO70} (see also \cite{Saito98}).   Consider
the rational function
$$\Phi(T):=T^{-d} 
\frac{(T^d-T^{q_1})(T^d-T^{q_2})(T^d-T^{q_3})}{(T^{q_1}-1)(T^{q_2}-1)(T^{q_3}-1)} .$$
By \cite[(1.3) Theorem]{Saito87} there exist finitely many
integers
$m_1, \ldots , m_\mu$ such that
$$\Phi(T)= T^{m_1} + \ldots + T^{m_\mu}.$$
Then
$\omega_i := \exp (2 \pi \sqrt{-1} m_i/d)$, $i=1, \ldots, \mu$, are the zeros of
$\phi_M(t)$. 
Put $\Lambda_k:= \omega_1^k + \ldots + \omega_\mu^k$
for $k \in \NN$ and let $\omega := \exp (2 \pi \sqrt{-1}/d)$.
Then one has $\Lambda_k=\Phi(\omega^k)$. From this one can derive that
$$\Lambda_k=\left( \delta(kq_1 \mbox{ mod }d)\frac{d}{q_1}-1 \right)
\left( \delta(kq_2 \mbox{ mod } d)\frac{d}{q_2}-1 \right)
\left( \delta(kq_3 \mbox{ mod } d)\frac{d}{q_3}-1 \right)$$
where $\delta$ is the delta function, i.e., $\delta(0):=1$ and $\delta(x):=0$ for $x
\neq 0$. The numbers $\Lambda_k$ and $\chi_m$ are related by the formula
$$\Lambda_k= \sum_{m|k} m \chi_m.$$

\noindent {\em Proof of Theorem~\ref{phi_M3}.} 
We have
$$\tilde{\phi}_A(t) = \frac{(1-t^d)(1-t^{\alpha_1}) \cdots (1-t^{\alpha_r})}
{(1-t)^{2g-2+r}(1-t^{q_1})(1-t^{q_2})(1-t^{q_3})}. $$ 
From Proposition~\ref{OW}
we conclude that $q_i|d$ or $q_i=\alpha_j$ for some $j$,
$1 \leq j \leq r$, and that $\alpha_i|d$ or $\alpha_i=q_j$ for some $j$, $1 \leq j \leq
3$. Therefore we may assume that
$$\tilde{\phi}_A(t) = \frac{(1-t^d)
\prod_{\alpha_i|d}(1-t^{\alpha_i})}{(1-t)^{2g-2+r}\prod_{q_j|d}(1-t^{q_j})}.$$
Hence
$$\tilde{\phi}_A^\ast(t)=
\frac{(1-t^d)^{2g-2+r}\prod_{q_j|d}(1-t^{d/q_j})}{(1-t)\prod_{\alpha_i|d}(1-t^{d/\alpha_i})}.$$
Denote by $\tilde{\Lambda}_k$ the sum of the $k$-th powers of the roots of
$\tilde{\phi}_A^\ast(t)=0$. For the proof of Theorem~\ref{phi_M3} we must show
that $\tilde{\Lambda}_k = \Lambda_k$ for all $k \in \NN$. 

(a) We first have
$$\tilde{\Lambda}_1= -1= \Lambda_1.$$

(b) Now suppose that $kq_{i_1} \equiv 0 \, \mbox{mod} \, d$  but $kq_{i_2} \not\equiv 0
\, \mbox{mod} \, d$, $kq_{i_3} \not\equiv 0 \, \mbox{mod} \, d$. Then we claim that
for all $i$ with $\alpha_i|d$ the number $\frac{d}{\alpha_i}$ does not divide
$k$. For suppose the contrary. By Proposition~\ref{OW} we have
$\alpha_i=(q_{j_1},q_{j_2})$ or $\alpha_i=q_{j_1}$. Now $\frac{d}{(q_{j_1},q_{j_2})}|k$
implies $d|kq_{j_1}$ and $d|kq_{j_2}$, which contradicts our assumption. But
$\alpha_i=q_{j_1}$ is only possible if $\alpha_i \! \not|d$.  Therefore we have shown
$$\tilde{\Lambda}_k = -1 + \frac{d}{q_{i_1}} = \Lambda_k.$$

(c) We now consider the case that $kq_{i_1} \equiv 0 \, \mbox{mod} \, d$, $kq_{i_2}
\equiv 0 \, \mbox{mod} \, d$, $kq_{i_3} \not\equiv 0 \, \mbox{mod} \, d$. By
Proposition~\ref{OW} and the same arguments as in (b), we see that $\alpha_i|d$
and
$\frac{d}{\alpha_i}|k$ only if $\alpha_i=(q_{i_1},q_{i_2})$. By
Proposition~\ref{OW} we therefore get
\begin{eqnarray*}
\tilde{\Lambda}_k & = & -1 +\delta(d \, \mbox{mod}\, q_{i_1})\frac{d}{q_{i_1}} +
\delta(d \, \mbox{mod}\, q_{i_2})\frac{d}{q_{i_2}} -   \sum_{\alpha_j=(q_{i_1},q_{i_2})}
\frac{d}{\alpha_j}\\
  & = & -1 + \frac{d}{q_{i_1}} + \frac{d}{q_{i_2}} - \frac{d^2}{q_{i_1}q_{i_2}} \\
  & = & -\left( \frac{d}{q_{i_1}}-1 \right) \left(\frac{d}{q_{i_2}}-1 \right)\\
  & = & \Lambda_k.
\end{eqnarray*}

(d) Finally, assume that $kq_1 \equiv 0 \, \mbox{mod} \, d$, $kq_2
\equiv 0 \, \mbox{mod} \, d$, and $kq_3 \equiv 0 \, \mbox{mod} \, d$. Since the
greatest common divisor of $q_1$, $q_2$, $q_3$ is 1, it follows that $d|k$. Then we have
\begin{eqnarray*}
\tilde{\Lambda}_k & = & (2g-2+r)d + \sum_{j=1}^3 \delta(d \, \mbox{mod} \, q_j)
\frac{d}{q_j} - \sum_{i=1}^r \delta(d \, \mbox{mod} \, \alpha_i) \frac{d}{\alpha_i} -1 \\
  & = & (2g-2+r)d + \sum_{j=1}^3 \frac{d}{q_j} - \sum_{i=1}^r \frac{d}{\alpha_i} -1.
\end{eqnarray*}
By Proposition~\ref{Wagreich3}(b) and (e) we get
\begin{eqnarray*}
\tilde{\Lambda}_k & = & \frac{Rd^2}{q_1q_2q_3} + \sum_{j=1}^3 \frac{d}{q_j} -1 \\
  & = &  \frac{d^2(d-q_1-q_2-q_3)}{q_1q_2q_3} + \sum_{j=1}^3 \frac{d}{q_j} -1 \\
  & = &  \left(\frac{d}{q_1}-1 \right) \left(\frac{d}{q_2}-1 \right)
\left(\frac{d}{q_3}-1 \right) \\
  & = & \Lambda_k.
\end{eqnarray*}
This completes the proof of Theorem~\ref{phi_M3}.

\section{The McKay correspondence} \label{KS}

In this section we shall derive Theorem~\ref{McKay} from the McKay correspondence.

Let $(X,x)$ be a Kleinian singularity. Then
${\cal D}=\PP^1(\CC)$, and
$\Gamma$ is a finite subgroup of $\mbox{Aut}(\PP^1(\CC)) = PGL(2,\CC)$. We may
assume that $\Gamma \subset PSU(2)
\cong SO(3)$. Up to conjugacy, there are five classes of such groups: (1) ${\cal
C}_{l+1}$, the cyclic group of order $l+1$, $l
\geq 1$, (2) ${\cal D}_{l-2}$, the dihedral group of order $2(l-2)$, $l \geq 4$, (3)
${\cal T}$, the tetrahedral group of order $12$, (4) ${\cal O}$, the octahedral group of
order
$24$, (5) ${\cal I}$, the icosahedral group of order $60$. It is well known that there
is a correspondence between these singularities and the irreducible root systems as
indicated in Table~\ref{table4.1}.

\begin{table}\centering
\caption{Kleinian singularities} \label{table4.1}
{\tabcolsep4pt
\begin{tabular}{|c|c|c|c|c|}
\hline
 & $\Gamma$ & $g; \alpha_1, \ldots , \alpha_r$ & Weights & $\pi_A$ \\
\hline
$A_{2n-1}$ & ${\cal C}_{2n}$ & $0; n, n$ & $1,n,n/2n$ & $2n/1$ \\
$A_{2n}$ & ${\cal C}_{2n+1}$ & $0; 2n \! + \! 1, 2n \! + \!1$ & $2,2n\! +
\!1,2n\! +\! 1/4n\! + \!2$ &
$4n\! + \! 2/2$ \\
$D_l$ & ${\cal D}_{l-2}$ & $0; 2, 2, l\! - \! 2$ & $2,l\! - \! 2,l\! - \!1/2(l \!-
\!1)$ & $2 \! \cdot   \! 2(l\! - \! 1)/1 \! \cdot \! (l \! - \! 1)$ \\
$E_6$ & ${\cal T}$ & $0; 2, 3, 3$ & $3,4,6/12$ & 
$2 \! \cdot \! 3 \! \cdot \! 12/1 \! \cdot \! 4 \! \cdot \! 6$ \\
$E_7$ & ${\cal O}$ & $0; 2, 3, 4$ & $4,6,9/18$ & 
$2 \! \cdot \! 3 \! \cdot \! 18/1 \! \cdot \! 6 \! \cdot \! 9$ \\
$E_8$ & ${\cal I}$ & $0; 2, 3, 5$ & $6,10,15/30$ & 
$2 \! \cdot \! 3 \! \cdot \! 5 \! \cdot \! 
30/1 \! \cdot \! 6 \! \cdot \! 10 \! \cdot \! 15$ \\
\hline 
\end{tabular}}
\end{table}

Let $\nu = -\frac{2}{R}$.
We have $A_k=S^{\nu k}(\CC^2)^\Gamma$ where $S^{\nu k}(\CC^2)$ denotes the $\nu k$-th symmetric
power of $\CC^2$ and the action of $\Gamma$ on $S^{\nu k}(\CC^2)$ is induced by its action
on $\PP^1(\CC)$.
We relate the series
$p_A(t)$ to another Poincar\'{e} series which is considered in \cite{Kostant85}. Let
$SU(2) \to PSU(2)$ be the usual double covering  and let $G
\subset SU(2)$ be the inverse image of $\Gamma \subset PSU(2)$. 
Let $\rho_m$ be the representation of $G$ on
$S^m(\CC^2)$ induced by its action on $\CC^2$. Let $\gamma_0, \ldots , \gamma_l$ be the
equivalence classes of irreducible finite dimensional complex representations of
$G$ where
$\gamma_0$ is the class of the trivial representation. For each integer $m \geq 0$ we
have a decomposition $\rho_m = \sum_{i=0}^l v_{mi} \gamma_i$ with $v_{mi} \in \ZZ$.
We associate to $\rho_m$ the vector $v_m=(v_{m0}, \ldots , v_{ml})^t \in \ZZ^{l+1}$.
As in \cite[p.~211]{Kostant85} we define
$$P_G(t):= \sum_{m=0}^\infty v_m t^m.$$
This is a formal power series with coefficients in $\ZZ^{l+1}$. We also put 
$P_G(t)_i:= \sum_{m=0}^\infty v_{mi} t^m$. 
Note that $v_{m0}$ is the dimension of the $G$-invariant subspace of
$S^m(\CC^2)$. 
If $-I \in G$, then we have $v_{m0}=0$ for $m$ odd. 
One has $-I \in
G$ if and only if $\Gamma$ is not a cyclic group of odd order which is equivalent to
$R=-1$.  Therefore we get
$$p_A(t^\nu)=P_G(t)_0.$$

J.~McKay \cite{McKay80} has observed
that if $\gamma : G \to SU(2)$ is the given 2-dimensional representation
of $G$ then the $(l+1)
\times (l+1)$-matrix $B=(b_{ij})$, defined by decomposing the tensor products $\gamma_j
\otimes \gamma = \bigoplus_i b_{ij} \gamma_j$ into irreducible components, satisfies
$B=2I-C$
where $C$ is the affine Cartan matrix of the corresponding root system.  Moreover, the
indexing is so that the additional vertex in the extended Coxeter-Dynkin
diagram corresponding to the matrix $C$ has index 0.

\addvspace{3mm} 

{\em Proof of Theorem~\ref{McKay}.}
From the Clebsch-Gordon formula one can derive that
$$Bv_m=v_{m+1}+v_{m-1}$$
for all non-negative integers $m$ where $v_{-1}=0$ \cite[(3.3.1)]{Kostant85}. This
can be reformulated as follows (cf.\ \cite[p.~222]{Kostant85}).
Let $V$ denote the set of all formal power series 
$x = \sum_{m=0}^\infty x_m t^m$
with $x_m \in \ZZ^{l+1}$. This is a free module of rank $l+1$ over the ring $R$ of
formal power series with integer coefficients. 
Then $x=P_G(t)$ is a solution of the following linear
equation in $V$:
$$((1+t^2)I-tB)x=v_0.$$
Let $M(t)$ be the matrix $((1+t^2)I-tB)$ and $M_0(t)$ be the matrix obtained by replacing
the first column of $M(t)$ by $v_0=(1,0, \ldots, 0)^t$. Then Cramer's rule yields
$$P_G(t)_0= \frac{\det M_0(t)}{\det M(t)}.$$
From \cite[Ch.~V, \S~6, Exercice 3]{Bourbaki68} we obtain in the case when $R=-1$ 
$$\det M(t) = \det (t^2I-c_a), \qquad \det M_0(t) =\det(t^2I-c)$$
where $c$ is the Coxeter element and $c_a$ is the affine Coxeter element of the
corresponding root system. In the case $\Gamma={\cal C}_{l+1}$, $l+1=2n$, we assume
that the numbering of $\gamma_0, \ldots , \gamma_l$ is so that
the vertices of the extended Coxeter-Dynkin diagram (which is a cycle)
corresponding to $\gamma_0, \ldots , \gamma_{n-1}$ are not connected with each other and
the same holds for the vertices corresponding to $\gamma_n, \ldots , \gamma_l$. 
Note that this differs from the numbering in
\cite{Bourbaki68} but agrees with the numbering used for the discussion of the affine
Coxeter element in the cases different from $A_l$ in \cite{Steinberg85}. (The case $A_l$
is excluded in that paper.) This proves Theorem~\ref{McKay}.

\addvspace{3mm}

For a polynomial 
$$\phi(t)= \prod_{m|h} (1 - t^m)^{\chi_m},$$
we use the symbolic
notation
$$\pi:= \prod_{m|h} m^{\chi_m}.$$
In the theory of finite groups, this symbol is known as a {\em Frame shape} (cf.\
\cite{CN79}). The Frame shapes $\pi_A$ corresponding to the polynomials $\phi_A(t)$ are
indicated in Table~\ref{table4.1}. 

\section{Generalizations to ICIS}

In this section we shall consider generalizations of Theorem~\ref{phi_M3} to certain ICIS.

Let $(X,x)$ be an ICIS in $\CC^4$ given by quasihomogeneous equations $g=0$ and $f=0$
of degrees $d_1$ and $d_2$ respectively. As above, let $M$ be
the Milnor lattice and $c: M \to M$ be the monodromy operator of $(X,x)$. 
By \cite{GH78} the characteristic polynomial of the monodromy
operator can be computed as above using an appropriate rational function $\Phi(T)$.
Similarly to the proof of Theorem~\ref{phi_M3} one can show:

\begin{theorem} \label{phi_M4a}
Let $(X,x)$ be a quasihomogeneous ICIS in $\CC^4$ with weights $q_1$, $q_2$, $q_3$, $q_4$
and  degrees $d_1$, $d_2$.  Assume that $g(z_1,z_2,z_3,z_4)=z_1z_4+z_2z_3$. Define  
$$\tilde{\phi}_A(t):=
\frac{\phi_A(t)(1-t^{d_2})}{(1-t)^{2g}(1-t^{d_1})}, \quad \phi_M^\flat(t):=
\frac{\phi_M(t)}{(1-t)}.$$ 
Then we have  $\tilde{\phi}_A^\ast(t)=\phi_M^\flat(t)$.
\end{theorem}

\begin{theorem} \label{phi_M4b}
Let $(X,x)$ be a quasihomogeneous ICIS in $\CC^4$ with weights $q_1$, $q_2$, $q_3$, $q_4$
and  degrees $d_1$, $d_2$.  Assume that either

(A) $g(z_1,z_2,z_3,z_4)=z_1^q+z_2z_3$
and $f(z_1,z_2,z_3,z_4)=f'(z_1,z_2,z_3)+z_4^p$ for some integers $p,q \geq 2$ where
$q|d_2$, or

(B) $g(z_1,z_2,z_3,z_4)=z_1^q+(z_2-z_3)z_4$
and $f(z_1,z_2,z_3,z_4)=az_1^q+z_2(z_3-z_4)$ for some $a \in \CC$, $a \neq 0,1$, and some
integer $q \geq 2$ and $p:=2$. 

Define  
\begin{eqnarray*}
\tilde{\phi}_A(t) & := &
\frac{\phi_A(t)(1-t^{d_2})^{p-1}(1-t^{\frac{d_1}{q}})(1-t^{\frac{d_2}{p}})}
{(1-t)^{2g}(1-t^{d_1})(1-t^{\frac{d_2}{q}})^p}, \\ 
\phi_M^\flat(t) & := & \frac{\phi_M(t)(1-t^q)^p}{(1-t)^{p-1}(1-t^{\langle p,q
\rangle})^{(p,q)}}. 
\end{eqnarray*}
Then we have
$\tilde{\phi}_A^\ast(t)=\phi_M^\flat(t)$.
\end{theorem}

Note that in the case $p=2$, $p|q$, the polynomial $\phi_M^\flat(t)$ in
Theorem~\ref{phi_M4b} reduces to the corresponding polynomial of Theorem~\ref{phi_M4a}.

In case (A) of Theorem~\ref{phi_M4b}, $(X,x)$ is a $p$-fold suspension  and we use the
following result which can be derived from \cite[Theorem~10]{ESt98}. Let $(X,0)$ be an
ICIS in $\CC^{n+2}$ of dimension $n$ given by a map germ $F=(g,f): (\CC^{n+2},0) \to
(\CC^2,0)$. Let $X'=g^{-1}(0)$ and assume that $(X',0)$ is an isolated singularity. Let
$p \in \NN$ , $p \geq 2$. The {\em $p$-fold suspension} of $(X,0)$ is the ICIS
$(\tilde{X},0)$ defined by $\tilde{F}=(\tilde{g}, \tilde{f}) : (\CC^{n+2} \times \CC,0)
\to (\CC^2,0)$ where $\tilde{g}(y,z)=g(y)$ and $\tilde{f}(y,z)=f(y)+z^p$ for $(y,z) \in
\CC^{n+2} \times \CC$. Let $\phi_M$, $\phi'_M$, and $\tilde{\phi}_M$ be the
characteristic polynomials of the monodromy operators of the singularities $(X,0)$,
$(X',0)$, and $(\tilde{X},0)$ respectively. Write
$$ \phi_M(t)=\prod_{m|h} (1-t^m)^{\chi_m}, \quad \phi'_M(t) = \prod_{k|h'}
(1-t^k)^{\chi'_k}.$$
By \cite[loc.cit.]{ESt98} we have
$$\tilde{\phi}_M(t)= \prod_{m|h} \frac{(1-t^{\langle m,p
\rangle})^{(m,p)\chi_m}}{(1-t^m)^{\chi_m}} \prod_{k|h'} (1-t^{\langle k,p
\rangle})^{(k,p)\chi'_k}.$$

In case (B), $(X,x)$ is a Brieskorn-Hamm ICIS. A Brieskorn-Hamm ICIS is a singularity 
$(V_B(w_1, \ldots , w_n),0)$ (for integers $w_1, \ldots , w_n$, $w_i
\geq 1$, $n \geq 3$) where
$$V_B(w_1, \ldots, w_n):=\{ z \in \CC^n \, | \, b_{i1}z_1^{w_1}+
\ldots + b_{in}z_n^{w_n}=0; \ i=1, \ldots , n-2 \}$$
and $B=(b_{ij})$ is a sufficiently general $(n-2) \times n$-matrix of complex
numbers. In our case $(X,x)$ is analytically isomorphic to the singularity
$(V_B(q,2,2,2),0)$ for a suitable matrix $B$.  The orbit invariants of a 
Brieskorn-Hamm ICIS  are given in
\cite[Theorem~2.1]{NR78} and the characteristic polynomial of the monodromy operator is
computed in \cite{Hamm72}. 

\addvspace{3mm}

\noindent {\bf Example 1} Let $(X,x)$ be a simply elliptic ICIS. Then $(X,x)$ is
one of the singularities indicated in Table~\ref{table4.2} (cf.\  \cite{Saito74}). Let
$q_1, \ldots , q_n$ 
be the weights and $d$ be the degree of the equation(s) of $(X,x)$. Put $w_i =
\frac{d}{q_{n+1-i}}$, $i=1, \ldots , n$. Then $(X,x) \cong (V_B(w_1, \ldots , w_n),0)$.
Hence we can apply Theorem~\ref{phi_M3}  or  
Theorem~\ref{phi_M4b} if $n=3$ or $n=4$ respectively. We obtain
$\psi_A(t)=(1-t)^2$, $\tilde{\phi}_A(t)=p_A(t)$, and $\phi_M(t)=p_A^\ast(t)$ if 
$n=3$ and $\phi_M^\flat(t)=p_A^\ast(t)$ otherwise. The polynomial
$\phi_A(t)=1+(b-2)t+t^2$ is the characteristic polynomial of a rotation of the Euclidian
plane by the angle $\frac{2\pi}{d}$ and hence of the Coxeter element corresponding to
the Coxeter graph \cite{Bourbaki68}
$$
\unitlength1cm
\begin{picture}(1.8,0.5) 
\put(0.1,0.1){\circle{0.2}}
\put(0.2,0.1){\line(1,0){1.3}}
\put(0.8,0.2){$d$}
\put(1.6,0.1){\circle{0.2}}
\put(1.8,0){.}
\end{picture} 
$$
\begin{table}\centering
\caption{Simply elliptic ICIS} \label{table4.2}
\begin{tabular}{|c|c|c|c|c|c|}
\hline
Name & $\{g;b;\}$ & Weights & $w_i$ & $\pi_M$ & $\pi_A$ \\
\hline
$\tilde{E}_8$ & $\{1; 1; \}$ & $1,2,3/6$ & $2,3,6$ & $ 2 \cdot 3 \cdot 6/1$ &
 $1 \cdot 6/2 \cdot 3$ \\
$\tilde{E}_7$ & $\{1; 2; \}$ & $1,1,2/4$ & $2,4,4$ & $2 \cdot 4^2/1$ & $4/2$ \\
$\tilde{E}_6$ & $\{1; 3; \}$ & $1,1,1/3$ & $3,3,3$ & $3^3/1$ & $3/1$ \\
$\tilde{D}_5$ & $\{1; 4; \}$ & $1,1,1,1/2,2$ & $2,2,2,2$ & $2^4/1$ & $2^2/1^2$ \\
\hline 
\end{tabular}
\end{table}

\end{document}